\documentclass[leqno, myheadings, twoside]{amsart}
\usepackage{amsmath, amsthm, amssymb, amscd, amsxtra,graphicx}
\usepackage{latexsym, amsfonts}
\usepackage{url}
\usepackage{texdraw}
\usepackage{epsfig}
\def\R{\mathbb{R}}

\def\1-1{(1.1)}

\makeatletter \@addtoreset{equation}{section}

\makeatletter \renewcommand{\@biblabel}[1]{#1.}

\theoremstyle{remark}

\begin{document}
\title [Inequalities and asymptotic formulas for eigenvalues]
{Some inequalities and asymptotic formulas for eigenvalues on
Riemannian  manifolds}
\author{Genqian Liu}


\maketitle Department of Mathematics, Beijing Institute of
Technology,
 Beijing 100081, People's Republic of China.
 \ \
E-mail address:  liugqz@bit.edu.cn

\vskip 0.46 true cm

\vskip 15 true cm

\begin{abstract} \  In this paper, we establish sharp inequalities
for four kinds of classical eigenvalues on a bounded domain of a
Riemannian manifold. We also establish asymptotic formulas for the
eigenvalues of the buckling and clamped plate problems.
 In addition, we give a negative
  answer to the Payne conjecture for the one-dimensional case.
   \end{abstract}

\vskip 1.39 true cm

\section{ Introduction}

\vskip 0.45 true cm

Let $(M,g)$ be an $n$-dimensional oriented Riemannian manifold, and
let $\Omega\subset M$ be a bounded domain with $C^{2}$-smooth
boundary $\partial \Omega$.
   We consider the following classical eigenvalue problems:
    \begin{eqnarray} \left\{ \begin{array}{ll} \triangle_g u+\mu u=0
      \quad \;\; &\mbox{in}\;\; \Omega,\\
  \frac{\partial u}{\partial \nu}=0\quad \;\;  &
   \mbox{on}\;\; \partial \Omega, \end{array} \right. \end{eqnarray}
     \begin{eqnarray} \label{1.1}  \left\{\begin{array}{ll}\triangle_g u+\lambda u=0
      \quad \;\; &\mbox{in}\;\; \Omega,\\
  u=0\quad \;\;  &  \mbox{on}\;\; \partial
  \Omega, \end{array} \right. \end{eqnarray}
  \begin{eqnarray}   \left\{ \begin{array}{ll} \triangle_g^2 u
  -\Gamma^2 u=0 \quad \;\; &\mbox{in}\;\; \Omega,\\
  u=\frac{\partial u}{\partial \nu}=0\quad \;\;  &  \mbox{on}\;\;
  \partial \Omega, \end{array} \right. \end{eqnarray}
  \begin{eqnarray} \left\{ \begin{array}{ll} \triangle_g^2 u+\Lambda
  \triangle_g u=0 \quad \;\; &\mbox{in}\;\; \Omega,\\
  u=\frac{\partial u}{\partial \nu}=0 \quad \;\;  &  \mbox{on}\;\;
   \partial \Omega, \end{array} \right. \end{eqnarray}
  where $\nu$ denotes the outward unit normal vector to $\partial
  \Omega$, and $\triangle_g$ is the Laplace-Beltrami operator
  that is given in local coordinates by the
  expression,
   \begin{eqnarray*} \triangle_g =\frac{1}{\sqrt{|g|}}\sum_{i,j=1}^n
   \frac{\partial}{\partial x_i} \left( \sqrt{|g|}\,g^{ij}
   \frac{\partial}{\partial x_j}\right).\end{eqnarray*}
Here $| g | : =  det(g_{ij})$ is the determinant of the metric
tensor, and $g^{ij}$ are the components of the inverse of the metric
tensor $g$.

\vskip 0.12 true cm

   (1.1) is the Neumann problem (see \cite{Ch${}_1$});
   (1.2) is the Dirichlet problem (see \cite{Ch${}_1$} or \cite{CH});
     (1.3) occurs in the treatment of the vibration problem for a
    clamped plate (see \cite{CH} and \cite{PS}), and (1.4) is the well-known buckling
    problem for a clamped plate
     (see \cite{Ch${}_1$}, \cite{CH}, \cite{Pa}, \cite{Pay} or \cite{PS}).
       In each of these cases, the spectrum is discrete
       and we arrange the eigenvalues in non-decreasing order
       (repeated according to multiplicity)
   \begin{eqnarray*}  & 0=\mu_1 <\mu_2 \le  \cdots \le
   \mu_k\le\cdots;\\
    & 0<\lambda_1 < \lambda_2\le  \cdots \le
   \lambda_k\le\cdots; \\
      & \Gamma_1^2 \le \Gamma_2^2 \le \cdots  \le \Gamma_k^2 \le\cdots; \\
  & \Lambda_1 \le\Lambda_2 \le \cdots \le \Lambda_k \le\cdots.\end{eqnarray*}
   The corresponding eigenfunctions are expressed as
      $v_1$, $v_2$, $v_3$, $\cdots$;
      $\;u_1$, $u_2$, $u_3$, $\cdots$;
      $\;U_1$, $U_2$, $U_3$, $\cdots$;
      and $W_1$, $W_2$, $W_3$, $\cdots$.

 For any bounded domain $\Omega$ in $(M,g)$, the variational formulation of
the Neumann and Dirichlet eigenvalue problems (in terms of Rayleigh
quotients, cf. Sect.VI.1 of \cite{CH}) immediately implies the
inequalities
\begin{eqnarray*}\mu_k\le \lambda_k, \quad \,
k=1,2,3,\cdots.\end{eqnarray*}
  Moreover,
P\'{o}lya \cite{Pol} proved  in 1952 that for $\Omega \subset
{\mathbb{R}}^2$
\begin{eqnarray}  \mu_2<\lambda_1,\end{eqnarray}
 answering a
question of Kornhauser and Stakgold \cite{KS}. In the case that
$\Omega$ is a bounded convex domain $\Omega\subset {\mathbb{R}}^2$
with a piecewise $C^2$-smooth boundary, Payne \cite{Pa} showed that
\begin{eqnarray} \mu_{k+2} < \lambda_k, \quad k=1,2,3,\cdots. \end{eqnarray}
   Levine and
 Weinberger \cite{LW}  proved that
 \begin{eqnarray} \mu_{k+n}<\lambda_{k}, \quad k=1,2,3,\cdots\end{eqnarray}
 for smooth bounded
 convex domains $\Omega\subset {\mathbb{R}}^n$ (cf. \cite{Av}), as well as
 \begin{eqnarray}\mu_{k+m}\le \lambda_k, \quad\,
 k=1,2,3,\cdots; \, 1\le m\le n\end{eqnarray}
 for arbitrary bounded
 convex domains. In 1991, Friedlander \cite{Fr} proved that
\begin{eqnarray}\mu_k\le \lambda_{k+1}, \; \; k=1,2,3,\cdots\end{eqnarray}
 when   $\Omega\subset {\mathbb{R}}^n$ is a bounded domain with a $C^1$-smooth boundary
$\partial \Omega$.  We also refer to Mazzeo \cite{Ma} for an
extension to certain smooth manifolds, and to Ashbaugh and Levine
\cite{AL} and Hsu and Wang \cite{HW} for the case of subdomains of
the $n$-sphere ${\mathbb{S}}^n$ with a smooth boundary and
nonnegative mean curvature. Finally, in 2004 Filonov  \cite{Fi}
proved
 strict inequality
\begin{eqnarray} \label{1.2} \mu_{k+1}<\lambda_k \;\,  (k=1,2,3,\cdots).
\end{eqnarray}
 when $\Omega\subset {\mathbb{R}}^n$ is a domain with finite volume, and with the
 embedding $W^{1,2}(\Omega)\hookrightarrow L^2(\Omega)$ compact.

\vskip 0.30 true cm

  In regard to the vibration problem of a clamped plate,
P\'{o}lya in \cite{Po} obtained that
   \begin{eqnarray*} \lambda_k\le
\Gamma_k, \quad k=1,2,3,\cdots\end{eqnarray*}
 for any bounded domain
in ${\mathbb{R}}^2$.
 This result had actually been improved to be
\begin{eqnarray} \lambda_k<\Gamma_k,  \quad
 k=1,2,3,\cdots \end{eqnarray}
 for bounded domain $\Omega\subset {\mathbb{R}}^2$ with smooth boundary
 by Weinstein \cite{W} (which is referred as
  Weinstein's inequality, see \cite{WS}).
 In \cite{AL}, Ashbaugh and Laugesen
 obtained the inequalities
 \begin{eqnarray} \lambda_1^2\le \lambda_1
 \lambda_2 \le \Lambda_1 \lambda_1 \le \Gamma_1
\le \Lambda_1^2.\end{eqnarray}

\vskip 0.30 true cm

  Concerning  the buckling problem of a clamped plate,
 in 1937 Weinstein \cite{W} proved the following strictly inequality
   \begin{eqnarray} \lambda_k<\Lambda_k, \quad \, k=1,2,3,\cdots
   \end{eqnarray} for bounded domain $\Omega \subset
   {\mathbb{R}}^2$ with smooth boundary.
   Payne \cite{Pa} in 1955 proved that for any bounded
    domain $\Omega\subset {\mathbb{R}}^2$ with $C^2$-smooth
    boundary
     \begin{eqnarray}\lambda_2\le \Lambda_1,\end{eqnarray}
  solving a conjecture of Weinstein.
 In \cite {Pa}, Payne further made the conjecture:
 \begin{eqnarray} \lambda_{k+1} \le \Lambda_k, \,
 \quad \mbox{for} \;\; k=1,2,3,\cdots\end{eqnarray}
  Note that this conjecture remains open in ${\mathbb{R}}^n$
  ($n\ge 2$); however, in this paper we
  give a negative answer to the Payne conjecture for the one-dimensional case
   (see Example 4.1).

  \vskip 0.45 true cm

 The first purpose  of this paper is to prove:

\vskip 0.30 true cm

\noindent  {\bf Theorem 1.1.} \ \  {\it  Let $(M,g)$ be an
$n$-dimensional oriented Riemannian manifold ($n\ge2$), and let
$\Omega\subset M$ be a bounded domain with $C^{2}$-smooth boundary.
Then
  \begin{eqnarray} \label{1.3} \mu_k< \lambda_k<\Gamma_k<\Lambda_k
  \quad \mbox{for}\;\; k=1,2,3, \cdots,\end{eqnarray}
 where $\mu_k$, $\lambda_k$,  $\Gamma_k^2$ and  $\Lambda_k$ are
  the $k$-th eigenvalues
 of the Neumann, Dirichlet,  clamped plate
 and  buckling problems for the domain $\Omega$, respectively.}

\vskip 0.35 true cm

 We also show by some examples that in the Riemannian manifold
setting, (\ref{1.3}) are the best possible inequalities for these
classical eigenvalue problems (see Remark 3.1).

\vskip 0.30 true cm

   H. Weyl in 1912,
 was the first to establish asymptotic formulas in ${\mathbb{R}}^n$
 for the Dirichlet and Neumann eigenvalues  (see \cite{We} or \cite{Wey}):
   \begin{eqnarray}  &\lambda_k \sim (2\pi)^2
   \left(\frac{k}{\omega_n |\Omega|}\right)^{2/n}, \quad \, k\to \infty,\\
  & \mu_k \sim (2\pi)^2 \left(\frac{k}{\omega_n
  |\Omega|}\right)^{2/n},
  \quad \;\, k\to \infty,\end{eqnarray}
 where $\omega_n$ is the volume of the unit ball in
 ${\mathbb{R}}^n$, and $|\Omega|$ is the $n$-dimensional Lebesgue
  measure of $\Omega\,$ ($\sim$ means the ratio of the RHS to the LHS approaches $1$ as $k\to \infty$).
 In the case of two-dimensional Euclidean space, Pleijel \cite{Pl} in 1950 given the
asymptotic formula for the
  eigenvalues of a clamped plate based on a
Carleman's method in \cite{Ca} and \cite{Car}.
   In 1967, Mckean and Singer \cite{MS} generalized Weyl's asymptotic
 formulas to a bounded domain of a Riemannian manifold by
 investigating asymptotic expansion of the trace of heat operator.

\vskip 0.13 true cm

 The second purpose of the paper is to establish
 the asymptotic formulas for the eigenvalues of the buckling and
 clamped plate problems on a Riemannian manifold. We have the following:

\vskip 0.39 true cm

\noindent  {\bf Theorem 1.2.} \ \  {\it  Let $(M, g)$ be
 an $n$-dimensional Riemannian manifold, and let
$\Omega\subset M$ be a
  bounded domain with $C^{2}$-smooth boundary.
   Then
  \begin{eqnarray} \label{1.5} \;\; \qquad \quad\,\Lambda_k  \sim
 (2\pi)^2 \left(\frac{k}{\omega_n (\mbox{vol}(\Omega))}\right)^{2/n},
 \;\quad  \; \mbox{as}\;\; k\to
  +\infty,\end{eqnarray}
  \begin{eqnarray} \label{1.6} \;\; \qquad \quad\,\Gamma_k  \sim
 (2\pi)^2 \left(\frac{k}{\omega_n (\mbox{vol}(\Omega))}\right)^{2/n},
 \;\quad  \; \mbox{as}\;\; k\to
  +\infty,\end{eqnarray}
  where $\mbox{vol}(\Omega)$ is the volume of $\Omega$.}

\vskip 0.30 true cm

  Weyl's asymptotic formulas and Theorem 1.2 show that for general
   bounded domain  with $C^{2}$-smooth boundary
   in a Riemannian manifold, the four kinds of classical quantities
 $\mu_k$, $\lambda_k$, $\Gamma_k$ and $\Lambda_k$ have the same
 asymptotic behavior as $k\to +\infty$.

\vskip 0.26 true cm

 The proof of Theorem 1.1 uses a key result (Lemma 2.1), which generalizes the
 classical Holmgren's uniqueness theorem (see \cite{Ra})
 to the Rienannian manifold, and a
  technique of \cite{Fi} by which  Filonov proved the inequalities
  (\ref{1.2}). In order to prove Theorem 1.2,
 we first consider the case of the Euclidean space and then obtain
 the version in Riemanniann manifold by applying metric
 expansion in normal coordinates system.
  The main method is to approximate $\Omega$
  by a union of subdomains whose boundary are piecewise smooth
   that has been suitably contracted.
  We thus get a lower
  estimate for the counting function of the buckling eigenvalues
  if these subdomains are open, disjoint and lie inside $\Omega$.
 On the other hand, an upper estimate had been given in \cite{MS}
 (see also p.441 of \cite{CH}) by
 investigating the Neumann and Dirichlet eigenvalue problems.
 Thus the desired result is proved.

 \vskip 1.48 true cm

\section{Some lemmas}

\vskip 0.49 true cm

The following several lemmas will be needed.

\vskip 0.29 true cm

\noindent  {\bf Lemma 2.1.} \ \  {\it Let $\Omega$ be a bounded
domain with $C^2$-smooth boundary in an $n$-dimensional Riemannian
 manifold $(M,g)$, and let $0\ne u\in W^{1,2}_0(\Omega) \cap W^{2,2}(
\Omega)$ be a solution of (1.2). Then $\frac{\partial u}{\partial
\nu}\big|_{\partial \Omega}$ does not vanish identically on
$\partial \Omega$.}

\vskip 0.36 true cm

\noindent {\it Proof.} \ \   Let $F(x, \xi)$ be a fundamental
solution for the
 Helmholtz  operator $\triangle_g +\lambda$ on $M$ (i.e., $F(x,\xi)$
 satisfying
 \begin{eqnarray} (\triangle_g + \lambda)F(x,\xi)= \delta_x(\xi), \end{eqnarray}
where $\triangle_g$ denotes the Laplace operator taken with respect
to the variables $\xi$, and $\delta_x (\xi)$ is the Dirac
$\delta$-function. More precisely, $(\triangle_g +\lambda)F(x, \xi)
=0$ with respect to $\xi\ne x$ for any fixed $x$).
 For $x\in M$, we choose the normal coordinates centered at $x$.
 Since $F(x,\xi)$ is singular at $\xi=x$ we cut out from $\Omega$ a
geodesic ball $B(x,\epsilon)$ contained in $\Omega$ with center $x$,
radius $\epsilon>0$.
  From $u\in W^{1,2}_0(\Omega)$, we find by the same argument as in
  Corollary 6.2.43 of \cite{Ha} that $u=0$ on $\partial \Omega$.
 Since $\triangle_g F(x, \xi)
 +\lambda  F(x, \xi)=0$ in $\Omega\setminus B(x, \epsilon)$,
 by Green's formula we obtain
 \begin{eqnarray*}  0 &=& \int_{\Omega\setminus B(x,\epsilon)}
  u (\xi) \big(\triangle_g F(x, \xi)   +\lambda  F(x,\xi)\big) dV_g(\xi) \\
  &=& \int_{\Omega \setminus B(x,\epsilon)} \big(\triangle_g
u (\xi)\big) F(x,\xi) dV_g(\xi)
   +\int_{\partial (\Omega\setminus B(x,\epsilon))} u (\xi) \frac{\partial
  F(x,\xi)}{\partial \nu_\xi} dS_g(\xi)   \\
 &&- \int_{\partial
(\Omega\setminus B(x,\epsilon))}  F(x, \xi)
 \frac{\partial   u (\xi) }{\partial \nu_\xi} dS_g(\xi)
  +  \lambda \int_{\Omega \setminus B(x, \epsilon)}   u (\xi) \,F(x,\xi)
  dV_g(\xi)
   \\  &= & \int_{\partial (\Omega\setminus B(x, \epsilon))}
 u(\xi) \frac{\partial F(x, \xi)}{\partial \nu_\xi} dS_g(\xi)
   - \int_{\partial (\Omega\setminus B(x,
 \epsilon))} F(x, \xi) \frac{\partial u (\xi)}{\partial \nu_\xi} dS_g(\xi)\\
 &=& -\int_{\partial B(x,\epsilon)}
  u (\xi) \frac{\partial  F(x, \xi)}{
\partial \nu_\xi} dS_g(\xi)
- \int_{\partial \Omega}
  F(x,\xi) \frac{\partial u(\xi)}{
\partial \nu_\xi} dS_g(\xi)  \\
 && + \int_{\partial B(x, \epsilon)}  F(x, \xi) \frac{\partial
 u (\xi)}{\partial \nu_\xi}dS_g(\xi),  \end{eqnarray*}
 where $dS_g(\xi)$ denotes the $(n-1)$-dimensional volume element, and
  $\frac{\partial}{\partial \nu_\xi}$ denotes the derivative in the direction of
 the outward unit normal vector $\nu_\xi$ at $\xi$.
  We now
wish to evaluate the limits  of the individual integrals in this
formula for $\epsilon\to 0$.
 On $\partial B(x, \epsilon)$, we have
$F(x,\xi) =F_1(\epsilon)+O(\epsilon)$ since we have used the normal
coordinates. From proof of Theorem 9.4 of \cite{Mc}, we know that
for $n\ge 2$,
 \begin{eqnarray} F_1(z)=F_0 (z) [1+f(z)]+h(z)
 \quad \, \text{as} \,\, |z|\to 0, \end{eqnarray}
where
 \begin{eqnarray} F_0 (z)=\left\{ \begin{array}{ll}
\frac{|z|^{2-n}}{n(2-n)\omega_n}
 & \, \mbox{for} \,\, n>2,\\
\frac{1}{2\pi}\log |z| &\, \mbox{for} \,\,
n=2,\end{array}\right.\end{eqnarray}
$$f(z)=O(|z|^2)$$  and
\begin{eqnarray} h(z)=\left\{ \begin{array} {ll} \mbox{const}
+O(|z|^2)  \, \, &  \mbox{for} \,\, n=2,\\
 0  \,\,  &   \mbox{for odd} \,\, n>2,\\
  \mbox{const} \times  \log (\sqrt{\lambda}|z|/2)[1+O(|z|^2)] \, \, &
 \mbox{for even} \,\, n>2; \end{array} \right. \end{eqnarray}
here $\omega_n$, as before, denotes  the volume  of the unit ball in
${\mathbb{R}}^n$, and the $O(|z|^2)$ terms are analytic functions of
$|z|^2$.

Under the normal coordinates, $$\xi=q(\epsilon,\eta)=\exp_x
\epsilon\eta,$$ $\epsilon\ge 0$, $\eta\in {\mathfrak{S}}_x=\{\eta
\in M_x\big||\eta|=1\}$, about $x$. As discussed in Section III.1 of
\cite{Ch${}_1$}, the volume element $dV_g$ is given by
$$ dV_g (q(\epsilon, \eta))=\sqrt{|g(\epsilon,\eta)|} \,d\epsilon \,d\mu_x (\eta),$$
where $d\mu_x$ is the standard $(n-1)$-measure on
${\mathfrak{S}}_x$; and  the $(n-1)$-dimensional volume element of
$\partial B(x, r)$ is given by
$$ dS_g(q(\epsilon, \eta)) =\sqrt{g(\epsilon, \eta)} d\mu_x (\eta).$$
  The discussion of Sections III.1 and XII.8 of \cite{Ch${}_1$} shows that
  $$\lim_{\epsilon\to 0} \frac{\sqrt{|g(\epsilon,\eta)|}}{\epsilon^{n-1}}=1.$$
  Since $u \in W^{2,2}(\Omega)$,
  by applying  local regularity of elliptic euqations (see,
  for example, Theorem A.2.1 of \cite{Jo}) repeatedly, we get
 that   $u\in W^{j,2}(\overline{B(x, \epsilon)})$ for all
 $j=1,2,3,\cdots$, which implies $u\in C^\infty(\overline{B(x,
 \epsilon)})$. It follows from (2.2) and
 (2.3) that for $\epsilon\to 0$,
\begin{eqnarray*} &&  \bigg|\int_{\partial
B(x,\epsilon)} F(x, \xi)\frac{\partial  u (\xi)}{\partial \nu_\xi}\,
dS_g(\xi) \bigg| \\
&& \quad \qquad \le \big(n\omega_n
\epsilon^{n-1}+o(\epsilon^{n-1})\big) |F_1(\epsilon) +O(\epsilon)|\;
\sup_{B(x,\epsilon)} |\nabla u| \to 0.\nonumber\end{eqnarray*}
 Furthermore,
\begin{eqnarray*}  \int_{\partial B(x, \epsilon)} & u (\xi)
  \frac{\partial F(x,\xi)}{\partial \nu_\xi} dS_g(\xi)  =
\big(\frac{\partial F_1(\epsilon)}{\partial \epsilon} +O(1)\big)\,
\int_{\partial B(x, \epsilon)}
 u (\xi) dS_g(\xi) \nonumber \\
  & \quad \, =\left(\frac{1}{n\omega_n \epsilon^{n-1}}+
 b(\epsilon)+O(1)\right) \int_{\partial B(x,\epsilon)} u
(\xi)dS_g(\xi)  \to u (x),\end{eqnarray*}
 where $b(\epsilon)$ satisfies $\lim_{\epsilon\to 0} \,
  n\omega_n \epsilon^{n-1} b(\epsilon)=0$.  Altogether,  we get
 \begin{eqnarray}  u(x)=-\int_{\partial \Omega} F(x,\xi)
 \frac{\partial u(\xi)}{\partial
\nu_\xi}dS_g(\xi). \end{eqnarray}
  Since $u$ does not vanish identically in $\Omega$,
  by the above formula we get that $\frac{\partial u}
  {\partial \nu_\xi}\big|_{\partial \Omega}$
    does not vanish identically on $\partial \Omega$. $\square$

\vskip 0.29 true cm

\noindent {\bf Remark 2.2.} \ \  (a) \ \  When $\Omega$ is a bounded
domain with $C^{2,\alpha}$-smooth boundary in a real analytic
Riemannian manifold $(M,g)$, Lemma 2.1 can  be immediately proved as
follows. Suppose by contradiction that $u\in W^{1,2}_0
 (\Omega)\cap W^{2,2}(\Omega)$ satisfies
  \begin{eqnarray*} \left\{ \begin{array}{ll} \triangle_g u+\lambda u=0
      \quad \;\; &\mbox{in}\;\; \Omega,\\
  u=0, \quad \frac{\partial u}{\partial \nu}=0
  \quad \;\;  &  \mbox{on}\;\; \partial \Omega. \end{array} \right. \end{eqnarray*}
  Since the elliptic operator $\triangle_g +\lambda$ has real analytic
 coefficients in local coordinates chart, it follows from Shauder's
 estimate (see, for example, Theorem 6.15 of \cite{GT}) that $u\in
 C^{2,\alpha} (\bar \Omega)$. Applying Holmgren's uniqueness theorem (see
  Theorem 2 of p.42 in \cite {Ra} or p.433 of \cite{Ta}), we obtain $u\equiv 0$ in $\Omega$.
  This contradicts the assumption that $u$
 does not vanish identically in $\Omega$.

  (b) \ \ When $M={\mathbb{R}}^n$, the proof of Lemma 2.1 is
  quite easy. Indeed, it follows
  form Rellich's formula for the Dirichlet eigenvalue (see
  \cite{Re}) that
  $$\lambda = \frac{ \int_{\partial \Omega}
  \sum_{m=1}^n \left(\frac{\partial u}{\partial \nu}
 \right)^2 x_m \nu_m  dS}{2\, \int_\Omega u^2 dx}, $$
   where $\nu(x)= (\nu_1(x), \cdots, \nu_n(x))$ with $x\in \partial
   \Omega$.
 Since $\lambda\ne 0$, we get that $\frac{\partial u}{\partial
 \nu}\big|_{\partial \Omega}$ cannot vanish identically on
 $\partial \Omega$.

 \vskip 0.33 true cm

\noindent  {\bf Lemma 2.3.} \ \  {\it  Let $(M,g)$ be an
$n$-dimensional Riemannian manifold ($n\ge 2$), and let
$\Omega\subset M$ be a bounded domain with $C^2$-smooth boundary.
Then,
  for any $\tau$ we have
 $$W^{2,2}_0(\Omega) \cap M_\tau= \{0\},$$
 where $M_\tau =\{ u\in W_0^{1,2} (\Omega)
 \cap W^{2,2}(\Omega)\big| \triangle_g u+\tau u =0 \;\, \mbox{in}
 \;\, \Omega\}$}.

\vskip 0.28 true cm

\noindent {\it Proof.} \ \ Let $v\in W^{2,2}_0 (\Omega)$.
 It follows from Corollary 6.2.43 of \cite{Ha} that
   \begin{eqnarray*}  \frac{\partial^j v}{\partial \nu^j}=0 \;\;
   \mbox{on}\,\, \partial \Omega,
   \;\mbox{for}\;\; 0\le j<\frac{3}{2}.\end{eqnarray*}
   Thus, for any $v\in W^{2,2}_0 (\Omega)\cap M_\tau$, we have
\begin{eqnarray*} \left\{ \begin{array}{ll} \triangle_g v+\tau v=0
\quad \;\; &\mbox{in}\;\; \Omega,\\
  v=0, \;\; \frac{\partial v}{\partial \nu}=0 \quad \;\;
   &  \mbox{on}\;\; \partial \Omega. \end{array} \right. \end{eqnarray*}
     By applying Lemma 2.1, we get  $v\equiv 0$ in $\Omega$. \;\;
     $\square$

\vskip 0.29 true cm

 Denote by $\sigma_N$ (respectively, $\sigma_D$,  $\sigma_P$ and
 $\sigma_B$)
the spectra of the Neumann (respectively, the Dirichlet, the clamped
plate and the buckling) problem for a bounded domain in Riemannian
manifold $(M,g)$.
  Let
 \begin{eqnarray*} & N^{(N)}(\tau) =\# \{\mu_k\in \sigma_N\big| \mu_k\le
  \tau\}, \quad \,  N^{(D)}(\tau) =\# \{\lambda_k\in \sigma_D\big| \lambda_k\le
  \tau\},\\
   &  N^{(P)} (\tau)=\# \{\Gamma_k^2\in \sigma_P\big| \Gamma_k\le
  \tau\}, \quad \,
   N^{(B)} (\tau)=\# \{\Lambda_k\in \sigma_B\big| \Lambda_k\le
  \tau\}
  \end{eqnarray*}
 be the counting functions of $\sigma_N$, $\sigma_D$,
 $\sigma_P$ and $\sigma_B$, respectively. Each eigenvalue is counted as many times as its
 multiplicity.

\vskip 0.34 true cm

 \noindent  {\bf Lemma 2.4.} \ \  {\it  For any
   $\tau$ we have
     \begin{eqnarray}  \label {2.1} \qquad &N^{(N)}(\tau)=\max\{\dim L\big| L\subset W^{1,2}
   (\Omega), \;  \int_\Omega |\nabla_g u|^2\, dV_g
     \le \tau \int_\Omega
   |u|^2 dV_g, \; u\in L\},\\
   \qquad \label{2.2} &  N^{(D)}(\tau)=\max\{ \dim L\big| L\subset W^{1,2}_0
   (\Omega), \; \int_\Omega |\nabla_g u|^2 dV_g \le \tau \int_\Omega
   |u|^2 dV_g, \; u\in L\},\\
    \label{2.3} &  N^{(B)} (\tau)=\max\{\dim L\big| L\subset W^{2,2}_0
   (\Omega), \, \int_\Omega |\triangle_g u|^2 dV_g \le \tau \int_\Omega
   |\nabla_g u|^2 dV_g, \, u\in L\}, \\
    \label{2.4}\quad \qquad &N^{(P)} (\tau)=\max\{\dim L\big| L\subset W^{2,2}_0
   (\Omega), \, \int_\Omega |\triangle_g u|^2 dV_g \le \tau^2 \int_\Omega
   |u|^2 dV_g, \; u\in L\},\end{eqnarray}
   where $\nabla_g u$ is the gradient of $u$ which has the
    expression in local coordinates
    $$\nabla_g u= \sum_{i,j=1}^n \big(g^{ij} \frac{\partial
    u}{\partial x_i}\big)\frac{\partial }{\partial x_j}.$$}

\vskip 0.29 true cm

 \noindent {\it Proof} \ \  (i) \ \ The argument proving (\ref{2.1}) and (\ref{2.2}) is completely
 analogous to the one used in the Euclidean space (see \cite{Gl} or \cite{ES}).
  Actually, the formulas (\ref{2.1}) and (\ref{2.2}) are known as Glazman's
    variational principle.

  (ii) \ \   For any fixed $\tau$, let $\Lambda_1, \cdots, \Lambda_k$
   be  all the buckling eigenvalues that are not greater than $\tau$.
   Then the corresponding buckling eigenfunctions $W_1, \cdots, W_k$ span a
   $k$-dimensional linear subspace $\Im_k$ of $W_0^{2,2}(\Omega)$
   (Suppose by contradiction that $W_m=c_1W_1 +\cdots +c_{m-1}W_{m-1} +c_{m+1}W_{m+1}+\cdots
   +c_kW_k$ for some $m$, where $c_1\cdots, c_{m-1}, c_{m+1},\cdots, c_k$ are constants.
   Therefore,  $\int_\Omega \nabla W_m\cdot
   (\nabla W_m-\sum_{i\ne m} c_i\nabla_g W_i)dV_g =0$.
   Noticing that \begin{eqnarray*} \int_\Omega \nabla_g W_i\cdot \nabla W_j \; dV_g =\left\{
  \begin{array}{ll} 1  \quad \;\; &\mbox{when}\;\; i=j,\\
   0   \quad \;\;  &  \mbox{when}\;\;
   i\ne j, \end{array} \right. \end{eqnarray*}
  we obtain $\int_\Omega |\nabla_g W_m|dV_g =0$, so that $W_m$ is a
  constant in $\Omega$. In view of $W_m\big|_{\partial \Omega}=0$ we get
   that $W_m\equiv 0$ in $\Omega$, which is a contradiction).
   It suffices to prove that the right-hand side of (\ref{2.3}) is also $k$.
  If it is not this case, then there exists a
  $(k+1)$-dimensional linear subspace $L_{k+1}$ of $W_0^{2,2} (\Omega)$ such that
   \begin{eqnarray} \int_\Omega |\triangle_g u|^2 dV_g \le \tau
  \int_\Omega |\nabla_g u|^2 dV_g \quad \; \mbox{for all} \;\; u\in
  L_{k+1}.\end{eqnarray}
  Thus, $E\cap L_{k+1}\ne 0$ for any linear subspace $E$ of $W_0^{2,2}(\Omega)$ with
  $codim(E)=k$.
  It follows from this and the variational formula
    \begin{eqnarray*} \Lambda_{k+1} = \sup_{E\subset W_0^{2,2}(\Omega), \; codim
    E=k} \bigg(\inf_{w\in E}\;\;\frac{\int_\Omega |\triangle_g w|^2 dV_g}
    {\int_\Omega |\nabla_g u|^2
    dV}\bigg)\end{eqnarray*}
    that $\Lambda_{k+1}\le \tau$, which is a contradiction. Therefore (\ref{2.3}) holds.

    (iii) \ \ The proof of (\ref{2.4}) is similar to (ii). $\;\;\square$

\vskip 0.36 true cm

 \noindent  {\bf Lemma 2.5.} \ \  {\it  Let $(M, g)$ be a Riemannian manifold, and let
   $\Omega\subset M$ be a bounded domain with $C^2$-smooth boundary. Suppose
   $\Omega_1, \cdots, \Omega_m$ are pairwise disjoint domains in $\Omega$,
   each of which has piecewise $C^{2}$-smooth boundary.
  Arrange all the buckling eigenvalues of $\Omega_1, \cdots, \Omega_m$
  in an increasing sequence
     \begin{eqnarray}   \Lambda^*_1 \le \Lambda^*_2 \le \Lambda^*_3 \le \cdots \end{eqnarray}
with each eigenvalue repeated according to its multiplicity, and let
\begin{eqnarray*}  \Lambda_1 \le \Lambda_2 \le \Lambda_3 \le \cdots\end{eqnarray*}
be the buckling eigenvalues for $\Omega$. Then we have for all
$k=1,2, 3, \cdots $
  \begin{eqnarray} \Lambda_k\le \Lambda^*_k.\end{eqnarray}  }

 \noindent {\it Proof} \ \
   For $j=1,2, \cdots, k$, let $\psi_j:\Omega \to {\mathbb{R}}^n$ be a
   buckling eigenfunction of $\Lambda^*_j$ when restricted to the appropriate
   subdomain, and identically zero, otherwise. Obviously,
   \begin{eqnarray*} \int_\Omega \nabla_g \psi_i \cdot \nabla_g \psi_j
   \, dV_g =\left\{ \begin{array}{ll} 1 \quad \, \mbox{if}\;\, i=j,\\
    0 \quad \, \mbox{if} \;\, i\ne j. \end{array}\right. \end{eqnarray*}
  Let $W_1, \cdots, W_{k-1}$ be the buckling eigenfunctions
  corresponding to eigenvalues $\Lambda_1,  \cdots,
  \Lambda_{k-1}$, respectively, which satisfy
\begin{eqnarray*} \int_\Omega (\nabla_g W_i \cdot \nabla_g W_j) dV_g
=\left\{ \begin{array}{ll} 1 \quad \, \mbox{if}\;\, i=j,\\
    0 \quad \, \mbox{if}\,\; i\ne j. \end{array}\right. \end{eqnarray*}
 Consider the functions $f$ of the form
\begin{eqnarray*}  f=\sum_{j=1}^k \beta_j \psi_j, \end{eqnarray*}
 where $f$ satisfies
 \begin{eqnarray} \label{2-5} \sum_{j=1}^k \beta_j \int_\Omega (\nabla_g W_i
 \cdot  \nabla_g \psi_j)\, dV_g =0, \quad \; i=1,2,\cdots, k-1. \end{eqnarray}
 If we think of $\beta_1, \cdots, \beta_k$ as unknowns and
 $\int_\Omega (\nabla_g W_i\cdot \nabla_g \psi_j) \, dV_g$ as given coefficients,
 then system has more unknowns than equations and a nontrivial
 solution of (\ref{2-5}) must exist.
 Applying Green's formula and the definition of $\psi_j$, we have
  \begin{eqnarray*} \int_\Omega (\triangle_g \psi_i)(\triangle_g \psi_j) dV_g
 &=& \int_\Omega  \psi_i(\triangle_g^2\psi_j) dV_g
= -\Lambda^*_j\int_\Omega \psi_i (\triangle_g \psi_j)dV_g \\&=&
\Lambda^*_j\int_\Omega (\nabla_g\psi_i\cdot \nabla_g \psi_j)dV_g
=\left\{ \begin{array} {ll} \Lambda^*_j \quad \; \mbox{if}\;\, i=j,\\
 0 \; \quad \; \; \mbox{if}\;\, i\ne j.\end{array}\right.\end{eqnarray*}
 Hence  \begin{eqnarray*}  \Lambda_k \int_\Omega |\nabla_g f|^2 dV_g \le
 \int_\Omega |\triangle_g f|^2 dV_g =\sum_{j=1}^k
  \Lambda^*_j \beta_j^2 \le \Lambda^*_k \int_\Omega |\nabla_g f|^2 dV_g,
  \end{eqnarray*}
 which implies the desired result. $\quad \;\; \square$

 \vskip 1.3true cm

\section{Inequalities of eigenvalues}

\vskip 0.45 true cm

\noindent {\it Proof of Theorem 1.1.} \ \ (i) \ \ We shall prove
$\mu_k<\lambda_k$ for all positive integer $k$. The proof is
 analogous to \cite{Fi}.  For any fixed $\tau$, it follows
from (\ref{2.2}) of Lemma 2.4 that there exists a subspace $F$ of
  $W^{1,2}_0(\Omega)$ such that dim$\,F=N^{(D)}(\tau)$ and
 \begin{eqnarray*}
 \int_\Omega |\nabla_g u|^2 dV_g \le \tau
\int_\Omega |u|^2 dV_g, \quad \; u\in F.\end{eqnarray*}
 Let $u\in F\cap M_\tau$, where $M_\tau=\{v\big|
 \triangle_g v+\tau v=0 \,\,\mbox{in}\,\, \Omega, \;\mbox{and}\, \,
 \frac{\partial v}{\partial \nu}=0 \,\, \mbox{on}\,\, \partial
 \Omega\}$. Since
$\partial \Omega\in C^{2}$, it follows from $u\in M_\tau$ and the
regularity of elliptic equations (see, for example, \cite{ADN}
  or Theorem 8.12 of \cite{GT}) that $u\in W^{2,2}(\Omega)$. From
$u\in W^{1,2}_0(\Omega)$, we get that $u=0$ on $\partial \Omega$, as
mentioned earlier. This implies that $u$ is also a Dirichlet
eigenfunction with eigenvalue $\tau$. By Lemma 2.1, we
 get that $\frac{\partial u}{\partial \nu}$ cannot vanish identically
in $\Omega$, which contradicts the fact that $\frac{\partial
u}{\partial \nu}=0$ on $\partial \Omega$.
 Thus $F\oplus M_\tau$ is a direct
 sum and we denote it by $G_\tau$.
 Let $u+v\in G_\tau\subset W^{1,2}(\Omega)$, where $u\in F$ and $v\in M_\tau$. We have
 \begin{eqnarray*}    \int_\Omega |\nabla_g (u+v)|^2 dV_g &=&\int_\Omega (|\nabla_g u|^2
  +|\nabla_g v|^2 +2 \nabla_g u\cdot \nabla_g v)dV_g\\
  &=&\int_\Omega (|\nabla_g u|^2
  +|\nabla_g v|^2 -2 u(\triangle_g v)dV_g\\
   &\le& \tau \int_\Omega |u+v|^2 dV_g,\end{eqnarray*}
 so that  \begin{eqnarray*} N^{(N)}(\tau) \ge \mbox{dim}\; G_\tau
 =N^{(D)}(\tau)
+\mbox{dim}\,M_\tau,\end{eqnarray*} Taking $\tau=\lambda_k$, we have
\begin{eqnarray*} \#\{\mu_j\in \sigma_N
   \big|\mu_j<\lambda_k\}=
   N^{(N)}(\lambda_k)-\mbox{dim}\, M_{\lambda_k} \ge N^{(D)}(\lambda_k)=k.\end{eqnarray*}
   That is, $\mu_k<\lambda_k$.

\vskip 0.3 true cm

 (ii) \ \ It follows from (\ref{2.4}) of Lemma 2.4
that there exists a subspace $H_\tau$ of $W_0^{2,2}(\Omega)$ such
that dim$\,H_\tau=N^{(P)}(\tau)$ and
\begin{eqnarray*} \int_\Omega |\triangle_g w|^2 dV_g \le \tau^2 \int_\Omega
|w|^2 dV_g, \quad \,\forall \, w\in H_\tau.\end{eqnarray*}
 Let $K_\tau=\{v\big|\triangle_g v
 +\tau v=0\,\, \mbox{in}\,\, \Omega, \,\mbox{and}\,\, v=0 \,\, \mbox{on}\,\,
 \partial \Omega\}$, and let $u\in H_\tau \cap K_\tau$.
  Since $u\in H_\tau \subset W^{2,2}_0(\Omega)$, we find by Corollary 6.2.43 of \cite{Ha} that
 $u=\frac{\partial u}{\partial \nu}=0$ on $\partial \Omega$.
    Lemma 2.1 implies that $u=0$ in $\Omega$, therefore,
    the sum $G_\tau:= H_\tau\bigoplus K_\tau$ is direct.
  Let $u=w+v\in G_\tau$, where $w\in H_\tau, \; v\in K_\tau$.
   It follows from Green's formula and Schwarz's inequality that for
   $w\ne 0$,
     \begin{eqnarray*} \bigg(\frac{\int_\Omega |\nabla_g w|^2
   dV_g}{\int_\Omega |w|^2 dV_g}\bigg)^2 \le \frac{\int_\Omega |\triangle_g w|^2
   dV_g}{\int_\Omega |w|^2 dV_g}. \end{eqnarray*}
 From this and the definition of $H_\tau$, we get
   \begin{eqnarray*} \int_\Omega |\nabla_g w|^2 dV_g \le \tau
   \int_\Omega |w|^2 dV_g.\end{eqnarray*}
  Note that
  \begin{eqnarray*}
    \int_\Omega |\nabla_g v|^2 dV_g =\tau \int_\Omega |v|^2 dV_g, \;\;
    \mbox{for}\;\; v\in K_\tau\end{eqnarray*}
     and
    $$ \int_\Omega \nabla_g
  w\cdot \nabla_g v\, dV_g=-\int_\Omega
  w(\triangle_g v)dV_g=\tau \int_\Omega
  w\,v\,dV_g.$$
Therefore, for any $u=w+v\in G_\tau \subset W_0^{1,2}(\Omega)$ we
have
  \begin{eqnarray*} \int_\Omega |\nabla_g(w+v)|^2 dV_g &=&
  \int_\Omega \big(|\nabla_g w|^2 +|\nabla_g v|^2 +2\nabla_g w\cdot\nabla_g
  v\big)dV_g\\
   &\le& \tau \int_\Omega
  |w+v|^2\,dV_g.\end{eqnarray*}
    For $0=w\in H_\tau$, there is equality in the above inequality.
  It follows that
   \begin{eqnarray*} N^{(D)} (\tau) \ge \mbox{dim}\;G_\tau =N^{(P)}
   (\tau) +\mbox{dim\;}K_\tau.\end{eqnarray*}
   By taking $\tau=\Gamma_k$, we obtain
   \begin{eqnarray*} \#\{\lambda_j\in \sigma_D
   \big|\lambda_j<\Gamma_k\}=
   N^{(D)}(\Gamma_k)-\mbox{dim}\;K_{\Gamma_k}\ge N^{(P)}(\Gamma_k)=k,\end{eqnarray*}
  Hence $\lambda_k<\Gamma_k$.

\vskip 0.3 true cm

(iii) \ \ For fixed $\tau>0$, (\ref{2.3}) of
   Lemma 2.4 implies that there exists a
 subspace $L_\tau$ of $W_0^{2,2}(\Omega)$ such that
 dim$\,L_\tau=N^{(B)}(\tau)$ and
 \begin{eqnarray*} \int_\Omega |\triangle_g w|^2 dV_g \le \tau \int_\Omega
 |\nabla_g w|^2 dV_g, \,\quad w\in L_\tau.\end{eqnarray*}
 Denote $J_\tau=\{v\big| \triangle_g^2 v-\tau^2 v=0
 \,\,\mbox{in}\,\, \Omega, \,\, \mbox{and}\, \,
  v=\frac{\partial v}{\partial \nu}=0\,\, \mbox{on}\,\, \partial \Omega\},$
  and put $G_\tau=L_\tau +J_\tau$. We shall prove  $L_\tau \cap J_\tau=\{0\}$.
  Suppose that $0\ne u\in L_\tau \cap J_\tau$.
  Then, in view of $u=0$ on $\partial \Omega$ we get
   that $\nabla_g u$ and $\triangle_g u$ don't vanish identically
   in $\Omega$.
  It follows from Green's formula and Schwarz's inequality that for
  any $u\in W^{2,2}_0(\Omega)$,
   \begin{eqnarray} \,\quad\;\; \quad \;\; \int_\Omega |\nabla_g u|^2 dV_g =\bigg|-\int_\Omega
  u(\triangle_g u)dV_g\bigg| \le \left(\int_\Omega  |u|^2
  dV_g\right)^{1/2} \left(\int_\Omega |\triangle_g u|^2
  dV_g\right)^{1/2},\end{eqnarray}
   i.e.,
  \begin{eqnarray} \label{3-2} \frac{\int_\Omega |\triangle_g u|^2 dV_g}
  {\int_\Omega |u|^2 dV_g} \le \left(\frac{\int_\Omega |\triangle_g u|^2 dV_g}{
  \int_\Omega |\nabla_g u|^2 dV_g}\right)^2, \quad \forall u\in W^{2,2}_0(\Omega). \end{eqnarray}
  Note that
 \begin{eqnarray} \label{3-3} \frac{\int_\Omega |\triangle_g u|^2
  dV_g}{\int_\Omega |\nabla_g u|^2 dV_g}\le \tau, \quad \,
  \forall \; u\in L_\tau \end{eqnarray} and
 \begin{eqnarray} \label{3-4} \tau^2=\frac{\int_\Omega |\triangle_g u|^2 dV_g}{\int_\Omega |u|^2
   dV_g}, \quad \forall \; u\in J_\tau.\end{eqnarray}
   Therefore, Schwarz's inequality in (3.1) is an equality, which implies
   there exists a constant $\beta$ such that $\triangle_g u+\beta u=0$ in $\Omega$.
  Since $u=0$ on $\partial \Omega$, it follows that $\beta>0$ and
  $u$ is a Dirichlet eigenfunction. Thus, we find by
  $\frac{\partial u}{\partial \nu}\big|_{\partial \Omega}=0$
   and Lemma 2.1 that
     $u\equiv 0$ in $\Omega$. This shows that the sum $L_\tau \oplus
  J_\tau$ is direct (we still denote the direct sum  by $G_\tau$).
  For any  $u=w+v\in G_\tau\subset W_0^{2,2}(\Omega)$,
  where  $w\in L_\tau, \; v\in J_\tau$, we have
  \begin{eqnarray} \label{3-5} \int_\Omega |\triangle_g(w+v)|^2 dV_g &=&
   \int_\Omega \big[|\triangle_g w|^2 +|\triangle_g v|^2
  +2w (\triangle_g^2 v)\big]dV_g\\
 &=&\int_\Omega (|\triangle_g w|^2 +|\triangle_g v|^2 +2\tau^2 w
   v)dV_g.\nonumber\end{eqnarray}
  By (\ref{3-2})---(\ref{3-5}),
    we arrive at  \begin{eqnarray*}  \int_\Omega |\triangle_g
  (w+v)|^2 dV_g\le \tau^2 \int_\Omega
  |w+v|^2\,dV_g. \end{eqnarray*}
  This implies that
   \begin{eqnarray*} N^{(P)} (\tau) \ge \mbox{dim}\,G_\tau =N^{(B)}
   (\tau)+\mbox{dim\,}J_\tau,\end{eqnarray*}
i.e.,
   \begin{eqnarray*} N^{(P)} (\tau) -\mbox{dim}\,J_\tau\ge N^{(B)}
   (\tau).\end{eqnarray*}
  Setting $\tau=\Lambda_k$, we see that
     \begin{eqnarray*} \#\{\Gamma_j^2\in \sigma_P
   \big|\Gamma_j<\Lambda_k\}=
   N^{(P)}(\Lambda_k)-\mbox{dim\;}J_{\Lambda_k}\ge N^{(B)}(\Lambda_k)=k,\end{eqnarray*}
   that is, $\Gamma_k<\Lambda_k$. $\;\; \square$

\vskip 0.38 true cm

\noindent {\bf Remark 3.1.} \ \  (i) \ \
  For a bounded domain of a Riemannian manifold, Mazzeo \cite{Ma}
  had  showed that \begin{eqnarray} \mu_k\le \lambda_k, \quad \,
   k=1,2,3,\cdots.\end{eqnarray} (Actually,
  Nazzeo proved that inequalities $\mu_{k+1}\le \lambda_k$ when $M$ is a
  Riemannian symmetric space of noncompact type).
  Therefore, the strict inequalities \begin{eqnarray} \label{2.5} \mu_k<\lambda_k, \quad \,
  k=1,2,3,\cdots\end{eqnarray}
     is an improvement of Mazzeo's result in the general Riemannian manifold.
   The following example shows  inequalities (\ref{2.5}) cannot be improved such that
   $\mu_{k+1}\le \lambda_k$ holds for $k=1,2,3,\cdots$.
   In fact, for the
   spherical cap of radius $\delta>\frac{\pi}{2}$ on the sphere of
   radius $1$ in ${\mathbb{R}}^n$, one has $\mu_2 (\Omega)>\lambda_1(\Omega)$
   (see, Theorem 3 of p.44 in \cite {Ch${}_1$}). This fact was also pointed out by
  Mazzeo in \cite{Ma}. Therefore our strict inequalities (\ref{2.5}) are
   sharp.

(ii)  \ \  The inequalities \begin{eqnarray} \label
{2.6}\lambda_k<\Gamma_k, \quad \, \mbox{for}\;\; k=1,2,3,\cdots
\end{eqnarray}  are a generalization of Weinstein's inequality to
$n$-dimensional
 Rienannian manifold. Here our proof is completely different from
that of \cite{W}.
  The inequalities (\ref{2.6}) cannot be improved to be $\lambda_{k+1}\le
  \Gamma_k$ for $k=1,2,3,\cdots$. Indeed, let $\Omega$ be the unit disk
  $\{x\in {\mathbb{R}}^2\big||x|<1\}$. Denote by $J_m(r)$ the Bessel function
  of order $m$ and by $j^{(l)}_m$  its $l$-th positive zero. Then
  the Dirichlet eigenfunctions are
    \begin{eqnarray*} \phi_{m,l}=a_{m,l}J_{m} (\sqrt{\lambda_{m,l}}\;r)
    \left\{\begin{array}{ll}
   \cos m \theta,\\ \sin m\theta,\end{array}\right.
   \quad m=0, 1,2,\cdots; \, l=1,2,3,\cdots,\end{eqnarray*}
and the corresponding eigenvalues are $\lambda_{m,l}
=(j_m^{(l)})^2$. Thus $\lambda_1(\Omega)\approx (2.4048)^2$,
$\lambda_2(\Omega)=\lambda_3(\Omega)\approx(3.832)^2$. It follows
from p.231 of \cite{PS} that $\Gamma_1 (\Omega)\approx (3.1962)^2$
(where $3.1962\cdots$ is the first zero of $J_0(r) I_1(r) +J_1 (r)
I_0(r), \, r>0$, and $I_m(r)$ is the modified Bessel function of
order $m$). This means that $\lambda_2(\Omega)>\Gamma_1(\Omega)$.

(iii)  \ \ For $k=2,3,4,\cdots$, our inequalities
$\Gamma_k<\Lambda_k$ ($k=2,3,4,\cdots)$ are completely new  even for
the case $M={\mathbb{R}}^n$.
 It is also sharp since it cannot be improved such that $\Gamma_{k+1}\le
 \Lambda_k$ holds for $k=1,2,3, \cdots$.  In fact, let $\Omega=\{x\in
 {\mathbb{R}}^2\big||x|<1\}$. Then we claim that $\Gamma_2(\Omega)>\Lambda_1(\Omega)$.
 Suppose by contradiction that $\Gamma_2(\Omega)\le
 \Lambda_1(\Omega)$.
 It follows from Theorem 1.1 that $\lambda_2(\Omega)
 <\Gamma_2(\Omega)$. Thus we get
 $\lambda_2(\Omega)<\Lambda_1(\Omega)$. However, for the unit disk
 $\Omega$, it must be $\lambda_2(\Omega)=\Lambda_1(\Omega)\approx(3.832)^2$. This is
 a contradiction, and the claim is verified.

 \vskip 1.3true cm

\section{Asymptotic formula for the buckling
 eigenvalues in ${\mathbb{R}}^n$}

\vskip 0.48 true cm

First, we consider the one-dimensional buckling problem:
\begin{eqnarray}  u''''(x) + \Lambda u''(x) =0, \quad \;  0\le x\le L,\\
   u(0)=u(L)=0, \; u'(0)=u'(L)=0.\end{eqnarray}
   It is easy to check that the general solution of (4.1) is
   $$u(x)= C_1+ C_2 x + C_3\cos \sqrt{\Lambda} \, x  +C_4\sin\sqrt{\Lambda} \,
   x.$$
 The boundary conditions yield the following equations for the
 coefficients:
     \begin{eqnarray*} \left\{ \begin{array}{ll} C_1=-C_3, \quad C_2=-\sqrt{\Lambda}\,  C_4,\\
   C_1+ C_2L +C_3\cos\sqrt{\Lambda}
 \, L + C_4 \sin\sqrt{\Lambda}\, L=0,\\
  C_2-\sqrt{\Lambda} \, C_3 \sin \sqrt{\Lambda}\, L +\sqrt{\Lambda} \, C_4
  \cos\sqrt{\Lambda}\, L=0.\end{array} \right.
  \end{eqnarray*}
In order that this system of equations has a nontrivial solution,
$\sqrt{\Lambda}$ must satisfy \begin{eqnarray*} \sin
\frac{\sqrt{\Lambda}\,L}{2} \left[ 2\sin \frac{\sqrt{\Lambda} \,
L}{2} - \sqrt{\Lambda}\, L \cos\frac{\sqrt{\Lambda}\,
L}{2}\right]=0.\end{eqnarray*}
   From the equation $\sin\frac{\sqrt{\Lambda}\, L}{2} =0$, we obtain that
   $$\Lambda_{1,k}(L) =\left(\frac{2k\pi}{L}\right)^2, \quad \, k=1,2,3,\cdots,$$
   and the associated eigenfunctions are
   $$u_{1,k}(L,x)= 1- \cos\frac{2k\pi}{L}\, x,   \quad \,  k=1,2,3,
   \cdots;$$
 From  the equation $2\sin \frac{\sqrt{\Lambda} \, L}{2} - \sqrt{\Lambda}\, L
\cos\frac{\sqrt{\Lambda}\, L}{2}=0$, we get that
   \begin{eqnarray} \label{4-3} \tan \frac{\sqrt{\Lambda}\, L}{2}=\frac{\sqrt{\Lambda}\,
   L}{2}.\end{eqnarray}
  If we denote by $\{\sqrt{\Lambda_{2,k}(L)}\big|k=1,2,3,\cdots\}$
   all the positive roots of
 (\ref{4-3}), then
 \begin{eqnarray*} u_{2,k}(L, x)&=& 1+\frac{\sqrt{\Lambda_{2,k}(L)} \sin
 (L \,\sqrt{\Lambda_{2,k}(L)})}{\cos (L\sqrt{\Lambda_{2,k}(L)})-1}\,x-\cos
  (\sqrt{\Lambda_{2,k}(L)}\,x) \\
    && -\,\frac{\sin (L\,\sqrt{\Lambda_{2,k}(L)})}
 {\cos (L\sqrt{\Lambda_{2,k}(L)}) -1}\sin (\sqrt{\Lambda_{2,k}(L)}\, x)\end{eqnarray*}
  is the buckling eigenfunction corresponding to the eigenvalue $\Lambda_{2,k}(L)$.
   By solving the system of equations
   \begin{eqnarray*} \left\{\begin{array}{ll}
    y=x \; & \\
    y=\tan x &,  \end{array} \right.
  \end{eqnarray*}
   we find that $\frac{2k\pi}{L}<\sqrt{\Lambda_{2,k}(L)}<\frac{(2k+1)\pi}{L}$ for all
   $k=1,2, 3,\cdots$.

\vskip 0.29 true cm

\noindent  {\bf Example 4.1. } From the above argument, we obtain
all the buckling eigenvalues for the interval $[0,L]$:
  \begin{eqnarray} \quad \;  \Lambda_1 =\left(\frac{ 2\pi}{L}\right)^2, \;\;
  \Lambda_2 =\Lambda_{2,1}(L), \quad  \Lambda_3 =\left(\frac{4\pi}{L}\right)^2, \;\;
  \Lambda_4=\Lambda_{2,2}(L), \cdots\cdots. \end{eqnarray}
  A simple calculation shows that the Dirichlet eigenvalues for the interval $[0,
  L]$ are
   \begin{eqnarray} \quad  \lambda_1 =\left(\frac{ \pi}{L}\right)^2, \;\;
   \lambda_2 =\left(\frac{2\pi}{L}\right)^2, \;\;
   \lambda_3 =\left(\frac{3\pi}{L}\right)^2, \;\;
   \lambda_4 =\left(\frac{4\pi}{L}\right)^2,
   \cdots\cdots, \end{eqnarray}
  and the corresponding Dirichlet eigenfunctions are $u_k(x)= \sin
  \left(\frac{k\pi x}{L}\right), \;  k=1,2,3,\cdots$.
   Recall that $\Lambda_{2,1}(L) <\left(\frac{3\pi}{L}\right)^2$, \  i.e.,
   $\Lambda_2 <\lambda_3$.  This shows that the Payne conjecture is not true
   for the one-dimensional case.

\vskip 0.30 true cm

 \noindent  {\bf Theorem 4.2.} \ \  {\it  Let
$\Omega$ be a
  bounded domain in ${\mathbb{R}}^n$ with $C^{2}$-smooth boundary.
  Then,
  \begin{eqnarray} N^{(B)} (\tau) =(2\pi)^{-n} \omega_n |\Omega| \tau^{n/2}
  \big(1+o(1)\big) \quad \;\mbox{as} \;\; \tau\to
  +\infty.\end{eqnarray}}

\noindent {\it Proof.} \ We give the proof for $n=2$ only,
indicating at its conclusion how it can be augmented to yield the
$n$-dimensional case.

  (i) \  Let $Q=[0,a]\times [0,b]$ be a rectangle in ${\mathbb{R}}^2$.
  Let us consider the buckling eigenvalue problem on $Q$:
\begin{eqnarray*}  \left\{ \begin{array}{ll}
 &\triangle^2 u + \Lambda \triangle u =0, \quad \;  \mbox{in}\;\; Q,\\
  & u=\frac{\partial u}{\partial \nu} =0,\quad  \; \quad \quad \;\mbox{on}\;\; \partial Q.
  \end{array} \right.\end{eqnarray*}
 It is easy to check that the buckling eigenvalues are
    \begin{eqnarray*} & \pi^2
   \left(\frac{(2l)^2}{a^2}+\frac{(2m)^2}{b^2}\right),\; \quad\,
   \frac{(2\pi l)^2}{a^2}+ \Lambda_{2,m}(b),\\
    & \Lambda_{2,l}(a)+ \frac{(2\pi m)^2}{b^2}, \,\quad\;
  \Lambda_{2,l}(a) +\Lambda_{2, m}(b),\end{eqnarray*}
  and the corresponding eigenfunctions---up to a
normalizing factor---are respectively
 \begin{eqnarray*} (1- \cos\frac{2l\pi}{a}\, x)(1 - \cos\frac{2m\pi}{b}\,y); \quad
   \quad   \end{eqnarray*}
\begin{eqnarray*} (1 - \cos\frac{2l\pi}{a}\, x)\left(1+\frac{\sqrt{\Lambda_{2,m}(b)} \sin
 (b\sqrt{\Lambda_{2,m}(b)})}{\cos (b\sqrt{\Lambda_{2,m}(b)})-1}\, y-\cos
  (\sqrt{\Lambda_{2,m}(b)}\, y) \right.\\
   \left. -\frac{\sin (b\sqrt{\Lambda_{2,m}(b)})}
  {\cos (b\sqrt{\Lambda_{2,m}(b)}) -1}\sin (\sqrt{\Lambda_{2,m}(b)} \,y)\right);
  \quad \qquad  \end{eqnarray*}
\begin{eqnarray*} & \left(1+\frac{\sqrt{\Lambda_{2,l}(a)} \sin
 (a\sqrt{\Lambda_{2,l}(a)})}{\cos (a\sqrt{\Lambda_{2,l}(a)})-1}\, x-\cos
  (\sqrt{\Lambda_{2,l}(a)}\, x)  -\frac{\sin (a\sqrt{\Lambda_{2,l}(a)})}
  {\cos (a\sqrt{\Lambda_{2,l}(a)}) -1}\sin (\sqrt{\Lambda_{2,l}(a)}
  \,x)\right)\\
  & \quad \qquad \; \times  (1 - \cos\frac{2m\pi}{b}\, y);
  \quad \qquad  \end{eqnarray*}
   \begin{eqnarray*} & \left(1+\frac{\sqrt{\Lambda_{2,l}(a)} \sin
  (a\sqrt{\Lambda_{2,l}(a)})}{\cos (a\sqrt{\Lambda_{2,l}(a)})-1}\; x-\cos
 (\sqrt{\Lambda_{2,l}(a)}\; x)
 -\frac{\sin (a\sqrt{\Lambda_{2,l}(a)})}{\cos (a\sqrt{\Lambda_{2,l}(a)})-1}
 \sin (\sqrt{\Lambda_{2,l}(a)} \; x)\right)
  \\
  &\times  \left(1+\frac{\sqrt{\Lambda_{2,m}(b)} \sin
  (b\sqrt{\Lambda_{2,m}(b)})}{\cos (b\sqrt{\Lambda_{2,m}(b)})-1}\; y-\cos
  (\sqrt{\Lambda_{2,m}(b)}\; y)
 -\frac{\sin (b \sqrt{\Lambda_{2,m}(b)})}{\cos (b \sqrt{\Lambda_{2,m}(b)})-1}
 \sin (\sqrt{\Lambda_{2,m}(b)} \; y)\right),\; \quad \;\\
   &\mbox{for all }\; \,l,m=1,2,3, \cdots.
     \end{eqnarray*}
    If the number of the buckling eigenvalues for $\Omega$ no greater than a bound $\tau$ is
 again denoted  by $N^{(B)}(\tau)$,
 then $N^{(B)}(\tau)=N^{(B)}_1(\tau)+N^{(B)}_2(\tau) +N^{(B)}_3(\tau)+N^{(B)}_4(\tau)$,
 where
    \begin{eqnarray*} & N^{(B)}_1(\tau)= \#\{(l,m)\big|
    \frac{(2\pi l)^2}{a^2} + \frac{(2\pi m)^2}{b^2}
    \le \tau, \quad  \; l=1, 2, \cdots; m=1, 2, \cdots\};\quad \;\,\\
   &N^{(B)}_2(\tau)= \#\{ (l,m)\big| \frac{(2\pi l)^2}{a^2}
    +\Lambda_{2,m}(b) \le \tau, \quad
   \;  l=1, 2, \cdots; m=1, 2, \cdots\};\quad\\
  &N^{(B)}_3(\tau)= \#\{ (l,m)\big| \Lambda_{2,l}(a) + \frac{(2\pi m)^2}{b^2}
  \le \tau, \quad
    \;  l=1, 2, \cdots; m=1, 2, \cdots\};\quad \\
     &\quad N^{(B)}_4(\tau)= \#\{(l,m)\big| \Lambda_{2,l}(a)
     + \Lambda_{2,m}(b)   \le \tau, \quad
    \;  l=1, 2, \cdots; m=1, 2, \cdots\}.\;\quad\, \;\end{eqnarray*}
    $N^{(B)}_1(\tau)$ is precisely the number of
    integral solutions of the inequality
    \begin{eqnarray*} \frac{(2l)^2}{a^2}
  +\frac{(2m)^2}{b^2} \le \frac{\tau}{\pi^2},\end{eqnarray*}
  here $l>0$, $m>0$.
  For sufficiently large $\tau$, the ratio of the area of this sector
  to the number of lattice points $\{(2l,2m)\}$ contained in it is
 arbitrarily close to $4$. Indeed,
 for each lattice point $(2l,2m)$, we associate with it a square of
 length 2 of each side such that the lattice point is the
 upper-right corner of the square. Then,
  the region consisting of these squares is
   contained in this sector of the ellipse;
   if we add more squares with the side length $2$ through which the
  ellipse passes (we denote the number by $R(\tau)$),
  then the region formed by all the squares contains the sector of the ellipse.
  Thus, we get
  \begin{eqnarray} 4 N^{(B)}_1(\tau) \le \tau \frac{ab}{4\pi} \le 4N^{(B)}_1(\tau)
  +4R(\tau).\end{eqnarray}
   As in p.431 of \cite{CH}, $R(\tau)-1$ is not greater than the arc length of the
   quarter ellipse, and this increases only with $(\max\{a,b\})\sqrt{\tau}\,$
   (Actually, the arc length of the quarter ellipse is
   $a\int_0^{\pi/2}\sqrt{1-\frac{a^2-b^2}{a^2}\cos^2 t}\, dt$ if
   $a>b$).
   Therefore, we have the asymptotic formula
   \begin{eqnarray*} \lim_{\tau\to +\infty} \frac{N^{(B)}_1(\tau)}{\tau}=
   \frac{ab}{16\pi}.\end{eqnarray*}
   More precisely, we can write
   \begin{eqnarray} \label{4-8}  N^{(B)}_1(\tau)=
   \frac{ab}{16\pi}\tau +\theta_1 c_1(\max\{a,b\})\sqrt{\tau},\end{eqnarray}
 where $c_1$ is a constant independent of $\tau$, and $-1<\theta_1<1$.

   Next, we shall estimate $N^{(B)}_2(\tau)$.
   Recall that $\big(\frac{2m\pi}{b}\big)^2<\Lambda_{2,m}(b)
   <\big(\frac{(2m+1)\pi}{b}\big)^2$ for all
   $m=1,2,3, \cdots$.
      Then, we have that
   \begin{eqnarray*}  N^{(B)}_2(\tau)\ge  \#\left\{(l,m)\big|
   \frac{(2l)^2}{a^2} +\frac{(2m+1)^2}{b^2}
   \le \frac{\tau}{\pi^2}, \;\;
     l=1, 2, \cdots; m=1, 2, \cdots\right\}. \quad
    \end{eqnarray*}
  Similar to the argument for $N^{(B)}_1(\tau)$,
    the ratio of the area of this sector
   to the number of lattice points  $\{(2l, 2(m+1))\}$ contained in it is
   arbitrarily close to $4$
     for sufficiently large $\tau$. In fact, if a square
   with the side length $2$ lying below and to the
   left of each lattice point is associated with it,
  then the region composed of these squares is contained in this sector of the
  ellipse.
  Thus, we get
  \begin{eqnarray*} 4N^{(B)}_2(\tau) +T(\tau)+ 4R(\tau)
  \ge \tau \frac{ab}{4\pi},\end{eqnarray*}
  where $R(\tau)$ is as before, and  $T(\tau)$ is the number of the unit squares (in this sector)
   whose bottom sides lie in
   the $x$-axis. Obviously, $T(\tau)\le a\sqrt{\tau}$.
   Therefore, we have the asymptotic relation
   \begin{eqnarray*} \lim_{\tau\to \infty} \frac{N^{(B)}_2(\tau)}{\tau}\ge
   \frac{ab}{16\pi}.\end{eqnarray*}
   In other words, we can write
   \begin{eqnarray} \label{4-9} N^{(B)}_2(\tau) \ge
   \frac{ab}{16\pi}\tau +\theta_2 c_2(\max\{a,b\})\sqrt{\tau},\end{eqnarray}
 where $c_2$ is a constant independent of $\tau$, and $-1<\theta_2<1$.
  With the similar argument as before, we can get
   \begin{eqnarray} \label{4-10} N^{(B)}_i(\tau)\ge
    \frac{ab}{16\pi}\tau +\theta_i c_i(\max\{a,b\})\sqrt{\tau},
    \quad \;\ \mbox{for}\;\; i=3,4.\end{eqnarray}
 It follows from (\ref{4-8})---(\ref{4-10}) that
 \begin{eqnarray} N^{(B)}(\tau)\ge
   \frac{ab}{4\pi}\tau +\theta c(\max\{a,b\})\sqrt{\tau}.\end{eqnarray}

 (ii) \   Let $\Omega$ be a domain which may be
decomposed into a finite number, say $h$, of squares $Q_1, Q_2,
\cdots, Q_h$ (or $n$-dimensional cubes in the case of $n$
independent variables) of side $a$. Such domains will be called
square-domains (or
 $n$-dimensional cube-domains). The area of $\Omega$ is then $|\Omega|=ha^2$
 (or its $n$-dimensional volume is
$|\Omega|=ha^n$). We consider the buckling problem for the domain
$\Omega$.
  From the previous argument, we have
  \begin{eqnarray} \label{4-12} N^{(B)}_{Q_j} (\tau) \ge \frac{a^2}{4\pi} \tau +
   \theta_j c_j a_j\sqrt{\tau}, \quad j=1,2, \cdots, h,\end{eqnarray}
where $N^{(B)}_{Q_j}(\tau)$ is the
   number of the buckling eigenvalues less than or equal
   to $\tau$ for the subsquare
   $Q_j$.
  It follows from Lemma 2.5 that for the square-domain, the $k$-th
 buckling eigenvalue $\Lambda_k$ for the domain $\Omega$ is at most
 equal to the $k$-th  number $\lambda_k^*$ in the sequence
 consisting of all the eigenvalues of the subdomains $Q_j$
(ordered according to increasing magnitude and taken with their
 respective multiplicity).
   Therefore we have
   \begin{eqnarray}  N^{(B)}(\tau)\ge N^{(B)}_{Q_1}(\tau)
   +\cdots+ N^{(B)}_{Q_h}(\tau), \end{eqnarray}
 where $N^{(B)}(\tau)$ is the number of the buckling eigenvalues less
than or equal to $\tau$ for the domain $\Omega$. Since the number
$N^{(B)}_{Q_j}(\tau)$ have the form given by inequalities
(\ref{4-12}), we get
\begin{eqnarray} N^{(B)}(\tau) \ge \frac{|\Omega|}{4\pi} \tau +\theta c |\partial \Omega|
\sqrt{\tau}, \end{eqnarray} where $-1<\theta<1$ and $c$ is a
constant independent of $\tau$.

If there are $n$ independent variables instead of two,  the
preceding discussion is still valid expect for the expressions
$N^{(B)}_{Q_j}(\tau)$.  It is easy to see that
 \begin{eqnarray*} N^{(B)}_{Q_j} (\tau)\ge \frac{1}{2^n} \frac{\omega_n
 a^n\tau^{n/2}}{\pi^{n}}
  +\theta_j c_j a^{n-1}\tau^{\frac{n-1}{2}}, \;\; j=1,2, \cdots, h.
\end{eqnarray*}
 where $\omega_n$ is the volume of the unit ball in ${\mathbb{R}}^n$.
 This implies that for an $n$-dimensional polyhedron of volume $|\Omega|$
consisting of a finite number $h$ of congruent cubes, we have
\begin{eqnarray} N^{(B)} (\tau)\ge \frac{\omega_n |\Omega|\tau^{n/2}}{(2\pi)^n}
 +\theta c |\partial \Omega|^{n-1}\tau^{\frac{n-1}{2}}.
\end{eqnarray}
 where
  $-1<\theta<1$ and $c$ is a constant
independent of $\tau$.

  (iii) \ We now consider the buckling problem for arbitrary bounded
domain with $C^2$-smooth boundary. With the above argument, it is
possible to obtain a lower bound for $N^{(B)}(\tau)$.

 Suppose the plane is
partitioned into squares of side $a$, inducing a decomposition of
the domain $\Omega$ into $h$ squares $Q_1, Q_2, \cdots, Q_h$ and $r$
boundary domains $G_1, G_2, \cdots, G_r$.   It follows from Lemma
2.5 that
  \begin{eqnarray}  \label{4-16}  N^{(B)}(\tau)\ge N^{(B)}_{Q_1}(\tau)+N^{(B)}_{Q_2}(\tau)
  +\cdots +N^{(B)}_{Q_h}(\tau);\end{eqnarray}
  furthermore, by (\ref{4-12}) we have
  \begin{eqnarray} \label{4-17} \quad  N^{(B)}_{Q_1}
  (\tau)+\cdots + N^{(B)}_{Q_h} (\tau) \ge
  \frac{ha^2}{4\pi} \tau +\theta c ha\sqrt{\tau} =
  \tau \left(\frac{ha^2}{4\pi} +\frac{\theta c
  ha}{\sqrt{\tau}}\right),\end{eqnarray}
 where, as before, $-1<\theta<1$ and $c$
is a constant independent of $a$ and $\tau$.

 By applying (\ref{1.3}) of Theorem
     1.1 we obtain that
     \begin{eqnarray} \label{4-18}  N^{(N)}(\tau)\ge N^{(D)}(\tau)\ge N^{(P)}(\tau) \ge
     N^{(B)}(\tau), \quad \; \mbox{for any } \tau.\end{eqnarray}

On the other hand, for the same partition of the domain $\Omega$, it
follows from p.441 of \cite{CH} that
 \begin{eqnarray}  \label{4-19} N^{(D)}(\tau) &\le& N^{(N)}_{Q_1}(\tau) + N^{(N)}_{Q_2}(\tau) +\cdots +
   N^{(N)}_{Q_h}(\tau) +N^{(N)}_{G_1} (\tau) +\cdots + N^{(N)}_{G_r}(\tau) \nonumber\\
  &\le&  \tau \left[ \left(\frac{ha^2}{4\pi} + \theta'_1 c'_2 r a^2\right)
  +(\theta'_2 c'_1 ha +\theta'_1c'_3
  ra)\frac{1}{\sqrt{\tau}}\right],\end{eqnarray}
 where $N^{(N)}_{Q_j}(\tau)$ and $N^{(N)}_{G_j}(\tau)$ are the numbers of the
  Neumann eigenvalues less than or equal to $\tau$ for $Q_j$ and $G_j$, respectively.
 From (\ref{4-16})---(\ref{4-19}), we obtain
  \begin{eqnarray}  \label{4-20}  &
 \tau \left[ \left(\frac{ha^2}{4\pi} + \theta'_1 c'_2 r a^2\right)
  +(\theta'_2 c'_1 ha +\theta'_1c'_3
  ra)\frac{1}{\sqrt{\tau}}\right] \qquad\qquad\qquad \\
   & \;\;\quad \quad \quad \;\;\ge    N^{(D)}(\tau)\ge N^{(P)}(\tau) \ge N^{(B)}(\tau)
  \ge \tau\left(\frac{ha^2}{4\pi} +\frac{\theta c h a}{\sqrt{\tau}}\right).\nonumber
     \end{eqnarray}
  Note that  $ar<\tilde{c}$. Hence, for sufficiently small $a$, we see that $a^2r$ and
  $\big|ha^2 -|\Omega|\big|$ are arbitrarily small. It follows from these inequalities that
\begin{eqnarray*}    \lim_{\tau\to \infty} \frac{4\pi
N^{(D)}(\tau)}{\tau |\Omega|}=\lim_{\tau\to \infty} \frac{4\pi
N^{(P)}(\tau)}{\tau |\Omega|}=\lim_{\tau\to \infty} \frac{4\pi
N^{(B)}(\tau)}{\tau |\Omega|}=1\end{eqnarray*} since we may take a
 sufficiently small fixed $a$ (the quantity $a$ can be arbitrarily
 chosen) such that the factor of $\tau$ in (\ref{4-20}) arbitrarily
close to the value $|\Omega|/4\pi$ for sufficiently large $\tau$.

   (iv) \ With the similar way as for the plane, we can get the
   desired result for the buckling eigenvalues in the $n$-dimensional case.
 $\;\;\square \;\;\; $

 \vskip 1.38 true cm

\section{Asymptotic formulas in Riemannian manifolds}

\vskip 0.46 true cm

 \noindent  {\bf Theorem 5.1.} \ \  {\it Let $(M, g)$ be
 an $n$-dimensional Riemannian manifold, and let
$\Omega\subset M$ be a
  bounded domain with $C^{2}$-smooth boundary.
  Then,
  \begin{eqnarray} \label{5.1} N^{(B)} (\tau) \sim
  (2\pi)^{-n} \omega_n (\mbox{vol}(\Omega)) \tau^{n/2}
   \quad \;\mbox{as} \;\; \tau\to
  +\infty,\end{eqnarray}
 \begin{eqnarray} \label{5.2} N^{(P)} (\tau) \sim (2\pi)^{-n} \omega_n
(\mbox{vol}(\Omega)) \tau^{n/2}
   \quad \;\mbox{as} \;\; \tau\to
  +\infty.\end{eqnarray} }

\vskip 0.30 true cm

 \noindent {\it Proof.} \  For any $x_0\in M$, we consider a geodesic,
  normal coordinates system at $x_0$.
 Under the normal coordinates one can expand the metric as follows:
\begin{eqnarray*} g_{ij}=\delta_{ij} -\frac{1}{3} \sum_{k,l=1}^n
R_{ikjl} x^k x^l +O(|x|^3)\end{eqnarray*} and
\begin{eqnarray*}\sqrt{\mbox{det}(g_{ij})} =1-\frac{1}{6} \sum_{i,j=1}^n R_{ij} x^i
x^j +O(|x|^3),\end{eqnarray*}
 where $R_{ikjl}$ and $R_{ij}$ are, respectively, the
 components of the curvature tensor and Ricci tensor associated with $g$;
 this is accomplished by
applying the exponential map to the tangent space at $0$ to obtain
coordinates on a patch and then fixing things up outside (see
\cite{Pe}, p.59 of \cite{CLN} or Chapter 10 of \cite{Cart}).
  We let $B_{x_0} (\varrho)$
  be the ball on which this coordinates system is defined.
  We can choose $\varrho$ sufficiently small such that
 in $B_{x_0} (\varrho)$, the eigenvalues of $g_{ij}$ and $g^{ij}$ are between
 $(1+\epsilon(\varrho))^{-1}$ and $(1+\epsilon(\varrho))$, and furthermore
   $dV_g =\sqrt{\mbox{det}(g_{ij})} \, dx$ where
 $\big(1+\epsilon(\varrho)\big)^{-1}< \sqrt{\mbox{det}(g_{ij})}
 < \big(1+\epsilon(\varrho)\big)$.
 Here $\epsilon(\varrho)$ is a positive function of variable $\varrho$, and
 $\epsilon(\varrho)\to 0$ as $\varrho\to 0$.
  Let $N$ be any compact sub-manifold in $M$ with
 \begin{eqnarray} \label{5-3}  Rc_g \ge -c \quad \; \mbox{on }N,\end{eqnarray}
  where $c$ is a positive constant.
  The classical Bochner-Lichnerowicz-Weitzenb\"{o}ck formula (see
  \cite{Lic}) reads
 \begin{eqnarray*} \int_N |\nabla^2_g u|^2dV_g =\int_N |\triangle_g
 u|^2 dV_g
- \int_N Rc_g (\nabla u, \nabla u) dV_g \quad \;\mbox{for any}\;\;
 u\in C^\infty_0(N),\end{eqnarray*}
  where $|\nabla^2_g u|^2$ is defined in an invariant ways as
  \begin{eqnarray*} |\nabla_g^2 u|^2 =\nabla^l \nabla^k u \,
  \nabla_l \nabla_k u = g^{pl} g^{kq} \left(\frac{\partial^2
  u}{\partial x^k \partial x^l} -\Gamma_{kl}^m \frac{\partial
  u}{\partial x^m}\right)\left(\frac{\partial^2
  u}{\partial x^p \partial x^q} -\Gamma_{pq}^r \frac{\partial
  u}{\partial x^r}\right)\end{eqnarray*}
 Together with (\ref{5-3}), it implies that
 \begin{eqnarray*} \int_N |\nabla^2_g u|^2dV_g \le
  \int_N |\triangle_g u|^2 dV_g +c\int_N |\nabla_g u|^2
 dV_g.\end{eqnarray*}
  Denote by ${\mathbb{B}}_0(\varrho)$ the ball of ${\mathbb{R}}^n$ with the center $0$ and
 radius $\varrho>0$, and denote  by $\triangle u$ and $\nabla u$  the
 usual the Laplacian and gradient of $u$ in ${\mathbb{R}}^n$.
    Passing in the coordinates system, we find by a similar way as in
p.135 of \cite{DHL} that
\begin{eqnarray*} \qquad (\triangle_g u)^2 \le (\triangle u)^2 +
  \tilde{\epsilon}(\varrho)|\nabla^2 u|^2 +\tilde{\epsilon}(\varrho)
  |\nabla u|^2, \quad \mbox{for}\;\;
  u\in C_0^2({\mathbb{B}}_0(\varrho))\end{eqnarray*}
  where $\tilde{\epsilon}(\varrho) \to 0$ as $\varrho\to 0$, while by
  the Bochner-Lichnerowicz-Weitzenb\"{o}ck formula,
  \begin{eqnarray*} \int_{{\mathbb{B}}_0(\varrho)} |\nabla^2 u|^2 dx =
  \int_{{\mathbb{B}}_0(\varrho)} (\triangle u)^2 dx.\end{eqnarray*}
      Note that for any $u\in C_0^2 ({\mathbb{B}}_0(\varrho))$,
   \begin{eqnarray*}  \int_{B_{x_0}(\varrho)} |\nabla_g u|^2 dV_g &=&
   \int_{{\mathbb{B}}_0(\varrho)} \sum_{i,j=1}^n g^{ij}
   \sqrt{\mbox{det}(g_{ij})} \,
\frac{\partial u}{\partial x_i}\, \frac{\partial u}{\partial x_j} dx
\\    &\ge&
  \int_{{\mathbb{B}}_0(\varrho)} (1+\epsilon(\varrho))^{-2} |\nabla u|^2
  dx. \nonumber
  \end{eqnarray*}
   Thus, we have that for any $u\in C_0^2 ({\mathbb{B}}_0(\varrho))$,
\begin{eqnarray} \label{5.4}\frac{\int_{B_{x_0} (\varrho)}
  (\triangle_g u)^2 dV_g }{\int_{B_{x_0} (\varrho)}
  |\nabla_g u|^2 dV_g} &\le&
    (1+\epsilon(\varrho))^2 (1+\tilde{\epsilon}(\varrho))\,
  \frac{\int_{
{\mathbb{B}}_0(\varrho)}
  (\triangle u)^2 dx}{\int_{{\mathbb{B}}_0(\varrho)}
    |\nabla u|^2 dx} \\
    &&  + (1+\epsilon(\varrho))^2
    \tilde{\epsilon}(\varrho).\nonumber
\end{eqnarray}
   We may always assume $\varrho$ is small enough
  such that $\lambda_1 ({\mathbb{B}}_0(\varrho))>1$,
   where $\lambda_1 ({\mathbb{B}}_0(\varrho))$ is the
   first Dirichlet eigenvalue for  ${\mathbb{B}}_0(\varrho)$.
   Since the geodesic open balls $\{B_{x_0}(\varrho) \big| x_0\in M\}$
   cover $\bar \Omega$, it follows from Lebesgue's lemma (see, for example,
 Theorem 6.27 of \cite{PM}) that there exists a
 constant $\gamma>0$ such that if any subdomain
 $G\subset \bar \Omega$ satisfies
   $diam (G)<\gamma$, then  $G$ must be contained in some $B_{x_0}(\varrho)$.
   Let us part the domain
   $\Omega$ into $h$ subdomains $G_1, G_2,
 \cdots, G_h$ with piecewise $C^{2}$-smooth boundaries such
  that $diam (G_j)<\gamma$, $1\le j\le h$.
 It follows from  Lemma 2.5 that the $k$-th
 buckling eigenvalue $\Lambda_k$ for the domain $\Omega$ is not greater than
 the $k$-th number $\Lambda_k^*$ in the sequence
 consisting of all the buckling eigenvalues of the subdomains $G_j$
 (arranged according to increasing magnitude and taken with their
 respective multiplicity).
 Thus, we have
 \begin{eqnarray}\label{5.5}  N^{(B)}(\tau) \ge N_{G_1}^{(B)}(\tau) +N_{G_2}^{(B)}(\tau) +\cdots
 +N_{G_h}^{(B)}(\tau),\end{eqnarray}  where $N^{(B)}(\tau)$ and $N^{(B)}_{G_j}(\tau)$ are
   the numbers of the buckling eigenvalues less than or equal to $\tau$ for $\Omega$ and
$G_j$, respectively.
  For each subdomain $G_j$, we take a point $p_j\in G_j$ such that
  $G'_j\subset {\mathbb{B}}_0 (\varrho)$,
 where $G'_j=\{x'\in {\mathbb{R}}^n \big| x'={\mbox{Exp}}_{p_j}^{-1} x, \; x\in
  G_j\}$. Therefore, under normal coordinates at $p_j$,
   the inequality (\ref{5.4}) holds for any $u\in W^{2,2}_0(G'_j)$.
   This implies that
  \begin{eqnarray} \;\quad \quad\qquad \; \label{5.6}  \Lambda_k(G_j)\le (1+\epsilon(\varrho))^2
  (1+\tilde{\epsilon}(\varrho))
  \Lambda_k(G'_j)+
   (1+\epsilon(\varrho))^2 \tilde{\epsilon}(\varrho),\quad \; k=1,2,3,\cdots.\end{eqnarray}
  By Theorem 1.1 and the Faber-Krahn inequality (see, for example,
   Theorem 2 of p.87 in \cite{Ch${}_1$}), we have
  \begin{eqnarray*} \Lambda_k (G'_j)\ge \Lambda_1 (G'_j) \ge \lambda_1(G'_j)\ge \lambda_1 (
 {\mathbb{B}}_0(\varrho))>1, \quad \; 1\le j \le h. \end{eqnarray*}
  It follows from this and (\ref{5.6}) that \begin{eqnarray*}\Lambda_k(G_j)\le
   (1+\epsilon(\varrho))^2 (1+2\tilde{\epsilon}(\varrho))
   \Lambda_k(G'_j),
  \quad \, j=1,2,\cdots,h,\quad \, k=1,2,3,\cdots, \end{eqnarray*}
   so that
     \begin{eqnarray} \label{5.7} N^{(B)}_{G_j} (\tau) \ge N^{(B)}_{G'_j}
     \left(\frac{\tau}{(1+\epsilon(\varrho))^2(1+2\tilde{\epsilon}(\varrho))}\right), \, \; \;
  j=1,2,\cdots, h.\end{eqnarray}
   By Theorem 4.2, we have that
   \begin{eqnarray} \label{5.8} && N^{(B)}_{G'_j}
  \left(\frac{\tau}{(1+\epsilon(\varrho))^2(1+2\tilde{\epsilon}(\varrho))}\right)
        =(2\pi)^{-n} \omega_n |G'_j| \\
          && \quad \quad \quad \times \left(
 \frac{\tau}{(1+\epsilon(\varrho))^2(1+2\tilde{\epsilon}(\varrho))}\right)^{n/2}
     (1+o(1)) \quad \; \mbox{as}\;\; \tau\to \infty.\nonumber\end{eqnarray}
     It follows from (\ref{5.5}), (\ref{5.7}) and (\ref{5.8}) that, as $\tau\to
     +\infty$,
        \begin{eqnarray*} \;\;\quad \; N^{(B)}(\tau) &\ge&
     (2\pi)^{-n} \omega_n \sum_{j=1}^h |G'_j |
     \left(\frac{\tau}{(1+\epsilon(\varrho))^2
  (1 +2\tilde{\epsilon}(\varrho))}\right)^{n/2}(1+o(1)).
 \end{eqnarray*}
    Recall that  $dV_g =\sqrt{\mbox{det}(g_{ij})} \, dx\,$
    with  $\,(1+\epsilon(\varrho))^{-1} <\sqrt{\mbox{det}(g_{ij})}
    <(1+\epsilon(\varrho))$. We have
  \begin{eqnarray*}\big(1+\epsilon(\varrho)\big)^{-1} |G'_j|
   < vol(G_j)< \big(1+\epsilon(\varrho)\big) |G'_j|,\end{eqnarray*}
   so that
 \begin{eqnarray*}   \sum_{j=1}^h |G'_j | &\ge & \big(1+\epsilon(\varrho)\big)^{-1}
  \sum_{j=1}^h (\mbox{vol} G_j)\\ &= &
\big(1+\epsilon(\varrho)\big)^{-1} \big(\mbox{vol}
 (\Omega)\big).\nonumber \end{eqnarray*}
   This implies that
\begin{eqnarray} \label{5-8} N^{(B)}(\tau) &\ge&
     (2\pi)^{-n} \omega_n
     \big(1+\epsilon(\varrho)\big)^{-1}(\mbox{vol}(\Omega)) \\
     &&  \times \left(\frac{\tau}{(1+\epsilon(\varrho))^2
     (1+2\tilde{\epsilon}(\varrho))}\right)^{n/2}(
1+o(1)) \quad \; \mbox{as}\;\; \tau \to \infty. \nonumber
\end{eqnarray} Hence
\begin{eqnarray} \label{5.9}\lim_{\tau\to \infty} \frac{N^{(B)}(\tau)}{
\tau^{n/2}}\ge (2\pi)^{-n} \omega_n (\mbox{vol}(\Omega)).
\end{eqnarray}
  For, we may choose the quantity
$\varrho$ arbitrarily, and by taking a sufficiently small fixed
$\varrho$, make the factor of $\tau^{n/2}$ in (\ref{5-8})
arbitrarily close to $(2\pi)^{-n} \omega_n (\mbox{vol}(\Omega))$ for
sufficiently large $\tau$.

  On the other hand, it follows from (6) of \cite{MS} that, for
  the bounded domain $\Omega$ in Riemannian manifold $(M,g)$,
  \begin{eqnarray} \label{5-11} \sum_{k=1}^\infty e^{-t\mu_k} &=&
  (4\pi t)^{-n/2}\left[\mbox{vol}(\Omega) +\frac{1}{4} \sqrt{4\pi t} \,
  (\mbox{vol}(\partial \Omega)) \right.\\
  && \left.+ \frac{t}{3} \int_\Omega R \,dV_g -
  \frac{t}{6} \int_{\partial \Omega} J \,dS_g +o(t^{3/2})\right],\nonumber\end{eqnarray}
 where $R$ is the scalar curvature at a point of $M$, and $J$ the
 mean curvature at a point of $\partial \Omega$. From (\ref{5-11}),
 we have \begin{eqnarray*} \int_{0-}^\infty e^{-t \tau} dN^{(N)}(\tau)
  =\sum_{k=1}^\infty e^{-t\mu_k} \sim (4\pi t)^{-n/2}
  (\mbox{vol}(\Omega)),\quad \mbox{as}\;\, t\to 0,\end{eqnarray*}
 where $N^{(N)}(\tau)$ is
   the number of the Neumann eigenvalues
   less than or equal to $\tau$ for $\Omega$,
   and $\int_{0-}^\infty e^{-t \tau} dN^{(N)}(\tau)$ is
   the Riemann-Stieltjes integral on $[0,
   +\infty)$
   (Note that $\int_{0-}^\infty e^{-t \tau} dN^{(N)}(\tau)$ means $\lim_{\delta\to 0+}
  \int_{-\delta}^\infty e^{-t \tau} dN^{(N)}(\tau)$).
 It follows from Proposition 3.2 of p.89 in \cite{Tal} that
  \begin{eqnarray*} N^{(N)} (\tau) \sim (2\pi)^{-n} \omega_n
  (\mbox{vol}(\Omega))\tau^{n/2}, \quad \, \mbox{as}\;\; \tau \to \infty,\end{eqnarray*}
  i.e.,
\begin{eqnarray} \label{5.10} \lim_{\tau\to \infty}
\frac{N^{(N)}(\tau)}{\tau^{n/2}} =(2\pi)^{-n} \omega_n
(\mbox{vol}(\Omega)),\end{eqnarray}
 By (1.16) of Theorem 1.1, we have that
 \begin{eqnarray} \label{5.11} N^{(N)}(\tau)\ge N^{(D)}(\tau)\ge N^{(P)}(\tau) \ge
     N^{(B)}(\tau), \quad \; \mbox{for any } \tau.\end{eqnarray}
 It follows from (\ref{5.9}), (\ref{5.10}) and (\ref{5.11}) that
 \begin{eqnarray}  \quad \;\; \lim_{\tau\to \infty}
\frac{N^{(B)}(\tau)}{\tau^{n/2}} &=& \lim_{\tau\to \infty}
\frac{N^{(P)}(\tau)}{\tau^{n/2}}  =  \lim_{\tau\to \infty}
\frac{N^{(D)}(\tau)}{\tau^{n/2}} \\ &=& \lim_{\tau\to \infty}
\frac{N^{(N)}(\tau)}{\tau^{n/2}}=(2\pi)^{-n} \omega_n
(\mbox{vol}(\Omega)).\nonumber    \quad \; \; \quad \square
\end{eqnarray}

\vskip 0.36 true cm

\noindent {\bf Remark 5.2.} \ \  (i) \ \  For the Dirichlet
 and Neumann eigenvalue problems, Seeley \cite{Se}  and Pham
 \cite{Ph} showed that if $\Omega\subset {\mathbb{R}}^n$ is a bounded domain
  with $C^\infty$-smooth boundary, then the following sharp remainder estimate holds:
 \begin{eqnarray*}  N^{(D)}(\tau) =  (2\pi)^{-n} \omega_n |\Omega|
 \tau^{n/2}(1+O(\tau^{-\frac{1}{2}})), \quad \, \mbox{as}\;\; \tau\to \infty.\end{eqnarray*}
  In \cite{Se}, Seeley used the method of hyperbolic equations
  which is the most precise of the known Tauberian methods.
  Seeley in \cite{See} has generalized the above result to $n$-dimensional Riemannian
 manifolds.

  (ii)  \ \  For the Dirichlet and Neumann eigenvalues of a
  bounded domain $\Omega$ in a smooth Riemannian manifold $M$,
  Ivrii (see \cite {Iv}) has established:
 \begin{eqnarray} \label{5-14} &\quad \quad \; N^\pm(\tau)=(2\pi)^{-n} \omega_n \cdot \mbox{vol}(\Omega)
 \cdot \tau^{n/2} \pm \frac{1}{4} (2\pi)^{-n+1} \omega_{n-1} \cdot
 \mbox{vol}(\partial \Omega)\cdot \tau^{(n-1)/2}\\
   &+o(\tau^{(n-1)/2}),
 \quad \, \mbox{as}\;\; \tau \to +\infty, \nonumber \end{eqnarray}
 under an additional assumption (roughly, that the set of ``multiply
 reflected periodic geodesics in $\bar \Omega$ is of measure
 zero''), where $N^{+}(\tau)$ and $N^{-}(\tau)$ denote the counting functions
of $\sigma_N$ and $\sigma_D$, respectively. Melrose \cite{Me}
 independently obtained the
 same asymptotic estimate (\ref{5-14}) for Riemannian manifolds with concave boundary.
  However,  Ivrii's method is no longer valid
  for the buckling and the clamped eigenvalues.

 \vskip 0.36 true cm

 \noindent {\it Proof of Theorem 1.2.} \ \  Taking  $\tau=\Lambda_k$
 in Theorem 5.1, we immediately obtain the
 conclusion of the theorem. \ \ $\square$

\vskip 2.6 true cm

\vskip 0.32 true cm

\end{document}